\UseRawInputEncoding
\pdfoutput=1 
\documentclass{article}

\usepackage[preprint]{main}
\usepackage{graphicx}
\usepackage{natbib}
\setcitestyle{numbers,square}

\usepackage[utf8]{inputenc} 
\usepackage[T1]{fontenc}    
\usepackage{hyperref}       
\usepackage{url}            
\usepackage{booktabs}       
\usepackage{amsfonts}       
\usepackage{nicefrac}       
\usepackage{microtype}      
\usepackage{xcolor}         
\usepackage{amsmath}

\title{Water Goes Where? A Water Resource Allocation Method Based on Multi-Objective Decision-Making}

\author{%
    Tongyue Shi\textsuperscript{\normalfont 1}$^{\ast}$ \quad
    Siyu Tao\textsuperscript{\normalfont 1}\thanks{Equal contribution.} \quad
    Haining Wang\textsuperscript{\normalfont 2}\thanks{Equal contribution.} \quad
    \\ 
    \textsuperscript{1}School of Computer Science and Technology, Soochow University, Suzhou, China\\
    \textsuperscript{2}School of Mathematical Sciences, Soochow University, Suzhou, China \\
    \texttt{ \{tyshi,sytao,hnwang\}@stu.suda.edu.cn} \\
}

\begin{document}

\maketitle

\begin{abstract}

For a long time, water and hydroelectric power are relatively important resources. Their rational distribution is closely related to regional agriculture, industry, residents, etc. In this paper, we mainly study the problem of allocation scheme for Glen Canyon Dam and Hoover Dam in the Colorado River Basin. Taking into consideration of various factors, we build models to achieve optimal scheduling.

Firstly, we propose the Water Strategy Decision Model, which can obtain different distribution methods for the different water levels. Also, we connect the two dams in series to consider the coupling effect between them and integrate this part into the main model. If there is no additional water supply, the water in the two lakes will be completely used up after $173$ days. When considering comprehensive water allocation, Mexico should also be taken into account, where the river has already stopped flowing in practice. At the same time, during peak water season, there will be $0.09$ MAF water that can flow into the Gulf of California from the Colorado River. In addition, we need to update the data every month, re-running the model parameters to accommodate the latest water allocation.

Secondly, we propose three criteria of allocation for water and power generation, namely, economic, social, and environmental benefits. The economic benefits mainly include the benefits of water used for industry, agriculture, and electricity, and we use a gray model to predict the benefit coefficients of water used respectively.  Social benefits mainly include the minimum shortage of water and electricity for agriculture, industry, and residents. Environmental benefits include mainly the number of pollutants in three aspects of wastewater discharges. For this multi-objective plan, the model is solved using a multi-objective ant colony genetic algorithm under the constraint of constant total water volume, and finally, the current reservoir capacities of the two lakes are input to derive the annual water supply to the five states.

Thirdly, based on the location and development characteristics of the five states, we obtain a water scheduling model based on the priority of geographic-industry characteristics. 

Fourthly, we can regard the model as a four-dimensional space of industrial, agricultural, residential and power generation water demand. The partial derivative calculation formula of multi-dimensional space can be used to obtain the results.

Finally,  we analyze the sensitivity of the model, and it shows that the model has strong adaptability and is easier to popularize. Moreover, we discuss the advantages and disadvantages of the models.

\end{abstract}

\section{Introduction} 
\subsection{Background} 

\subsubsection{Relevant Geographic Background}

\begin{itemize}
\item \textbf{the Colorado River Basin}
\end{itemize}

\textbf{The Colorado River}\cite{c3} is often referred to as the lifeblood of the American Southwest. The river originates in the Rocky Mountains of the western United States and then winds through seven states and approximately $1,400$ miles of stunningly diverse ecosystems before it reaches the below-sea-level desert expanses of Mexico. Parts of the seven states of Wyoming, Colorado, Utah, New Mexico, Arizona, California, and Nevada form the U.S. portion of the $243,000$ square-mile Colorado River Basin, with $2$\% of this area located internationally in Mexico’s Sonoran Desert.

\begin{figure}[h]\label{11}
\small
\centering
\includegraphics[width=10cm]{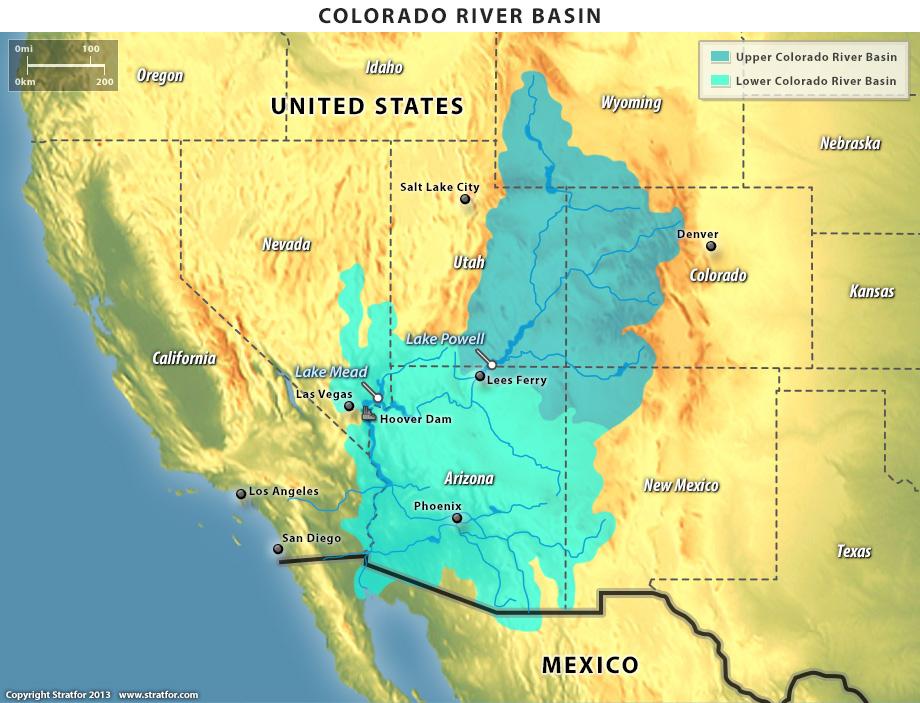}
\caption{Map of the colorado river basin }
\end{figure}

The 1922 Colorado River Compact\cite{c4}, created by these seven states, separated the basin into an upper and lower region with Lee’s Ferry just below Glen Canyon Dam as the point of division. Wyoming, Colorado, Utah, New Mexico, and the northern portion of Arizona make up the $109,800$ square-mile Upper Basin, while Arizona, California, and Nevada constitute the Lower Basin region as seen in Figure \ref{11}. The flow of it is not what it once was, as drought, over-apportionment, and ever-expanding urban development have depleted the supplies of this cherished resource to the point where it no longer reaches the Gulf of Mexico. Its future has become increasingly contentious and uncertain.

\begin{itemize}
\item \textbf{Glen Canyon Dam and Lake Powell}\cite{c1}
\end{itemize}

Glen Canyon Dam is a concrete arch-gravity dam on the Colorado River in northern Arizona, United States. It forms Lake Powell, one of the largest man-made reservoirs in the U.S. with a capacity of $27$ million acre-feet ($33 km^{3}$). The dam is also a major source of hydroelectricity, averaging over $4$ billion kilowatt-hours per year.

\begin{figure}[h]\label{12}
\small
\centering
\includegraphics[width=10cm]{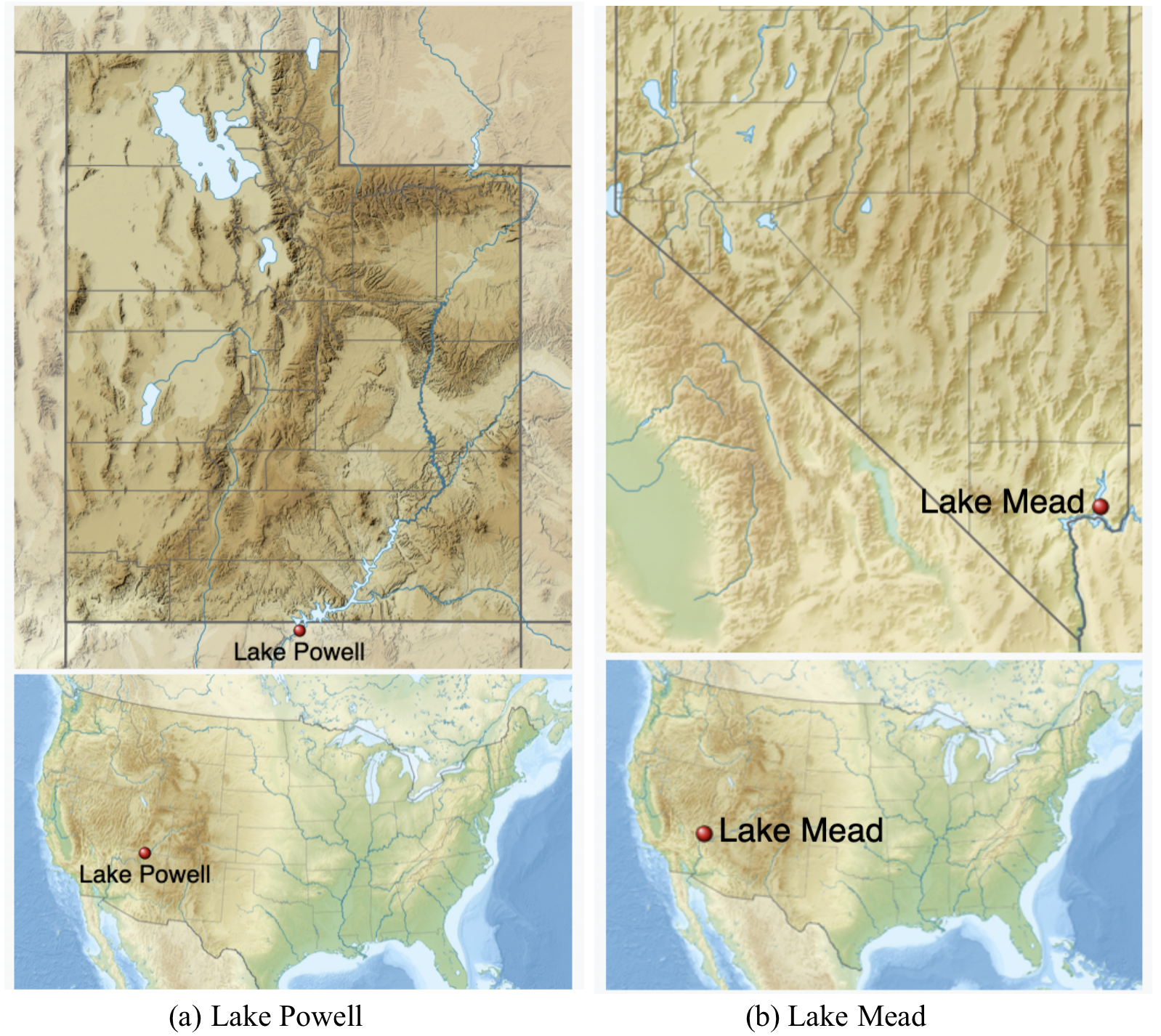}
\caption{The location of the two lakes and their state}
\end{figure}

\begin{itemize}
\item \textbf{Hoover Dam and Lake Mead}\cite{c2}
\end{itemize}

Hoover Dam is a concrete arch-gravity dam in the Black Canyon of the Colorado River, on the border between the U.S. states of Nevada and Arizona. Hoover Dam impounds Lake Mead, the largest reservoir in the United States by volume when full. The dam's generators provide power for public and private utilities in Nevada, Arizona, and California. 

\begin{table}[htb] 
    \centering
    \caption{Comparison of Glen Canyon Dam and Hoover Dam} 
    \begin{tabular}{ccccc} \hline
    \toprule
    Compare properties  &  Glen Canyon Dam & Hoover Dam  \\ 
\toprule

Height & 710 ft (220 m) &  726.4 ft (221.4 m) \\
Length & 1,560 ft (480 m) &   1,244 ft (379 m)\\
Elevation at crest & 3,715 ft (1,132 m) &  1,232 ft (376 m)\\
Reservoir & Lake Powell & 	Lake Mead \\
Total capacity & 27,000,000 af (33 $km^3$)&   28,537,000 af (35.200 $km^3$)\\
Catchment area & 108,335 sq.mi (280,590 $km^2$) & 167,800 sq.mi (435,000 $km^2$)\\
Surface elevation&  Current 3,554.01 ft (1,083 m)  &  1,219 ft (372 m)\\
Annual generation& 	4,717 GWh & 	3,300 GWh \\

\bottomrule \hline
\end{tabular}
\end{table}

\subsubsection{Problem Background}

Since the 19th century, the Colorado River Basin has been a developed agricultural area in the United States. More than $90$\% of the water consumption is used for agricultural irrigation. The irrigation scale in Arizona and California in the middle and lower reaches of the basin is particularly large. Among them, the Imperial Valley Irrigation District, an agricultural area in California, emerged in 1890, and the annual irrigation water volume is about $3.700$ billion $m^3$ from the Colorado River\cite{c5}. For the arid middle and lower reaches, the Colorado River is the lifeline for survival. Since the 19th century, when there were a lack of water conservancy projects and water rights systems, droughts, floods, and floods have frequently occurred in the basin, and water conflicts have existed for a long time. Although there are some contract clauses to arrange allocations, with problems such as high temperatures and scarce precipitation, water scarcity creates more and more problems. Therefore, how to allocate the relevant water resources and hydropower resources is the core issue to be studied in this thesis.

\subsection{Restatement of the Problem} 
Consider that we have to determine the best way to manage water usage and electricity production at \textbf{the Glen Canyon and Hoover dams} and propose a mathematical solution for the allocation of water in this region including the five U.S. states of \textbf{Arizona (AZ)}, \textbf{California (CA)}, \textbf{Wyoming (WY)}, \textbf{New Mexico (NM)}, and \textbf{Colorado (CO)}. We are supposed to address the following problems:

\begin{itemize} 
\item  \textbf{Problem 1}: Build a model that proposes a fixed set of water supply and demand conditions to analyze the operation of two dams with the water level of the two lakes as M and P. Determine the water volume that should be drawn from each lake to meet the stated demand and the time for the demand to be met without additional water supply. Calculate additional water volume which must be provided overtime to ensure that the fixed demand is met.

\item  \textbf{Problem 2}: Use the model from Problem 1 to propose the best way to address the competing interests of agriculture, industry, residential use, and electricity production and describe the criteria for addressing them. Propose the best allocation method and specify what criteria are used to address the interests.

\item  \textbf{Problem 3}: Build a model to solve problems when there is not enough water to meet the demand for water and electricity. 

\item  \textbf{Problem 4}: Explain the problem under the following conditions. 

\textbf{Condition1}: Agriculture, industry, and residential demand for water and electricity change over time. Figure out the trends when population, agriculture, and industry increase or decrease. 

\textbf{Condition2}:  Point out changes when the proportion of renewable energy exceeds the expected value. 

\textbf{Condition3}:  Describe changes when additional water and energy saving measures are implemented.

\item \textbf{Problem 5}: Write a one- to two-page article, representing the solution and giving your opinion, suitable for publication in \textit{Drought and Thirst}.

\end{itemize}

\subsection{Literature Review}
This area has always had relevant water rights treaties to allocate it. The water rights allocation treaty is an important subject of river basin management, but its implementation is inseparable from the protection and adjustment of engineering measures. Almost every water rights allocation will authorize a batch of water storage projects, water diversion projects and power generation projects. After the mega project is completed, its operation management and the environmental problems caused by it need to be solved by management means such as water rights distribution. The two complement each other and promote each other.

With the rapid development of water resources system analysis research and computer technology, optimization and simulation technology has made considerable progress in solving complex problems such as basin water resources planning, reservoir scheduling, water resources development and management.

Linear programming, nonlinear programming, dynamic programming and other optimization models have been widely used in water resources planning and management, and have played an important role in improving the level of water resources utilization and planning and scientific decision-making by decision-making departments.

In 1963, Dantzig studied the application of linear programming models to multi-objective problems\cite{c13}. In 1971, Joeres applied linear programming to Baltimore's multi-source water supply. In 1973, Dudley and Burt applied dynamic programming to the management of irrigation reservoirs, using the transition probability of the Markov chain to weight the recursive dynamic equation. In 1982, British scholars P. W Herberson and others, according to the characteristics of tidal power stations, considered the conflicting interests of multiple departments, and used a simulation model to simulate the distribution of fresh water in tidal bays, showing the superiority of simulation technology.

In addition, genetic algorithm and artificial neural network algorithm have been relatively mature in the field of water resources, and new intelligent algorithms such as ant colony algorithm, annealing algorithm, particle swarm algorithm have also been widely used in water resources basin. Applications.


\subsection{Review of Our Model and Work }

After data preprocessing, we first established a water resource scheduling model, considering the distance between the dam and the state, and using linear programming and numerical analysis to give the water resource scheduling solution.

Then, we use grey prediction and ant colony algorithm to solve the model parameters of water allocation for each industry by using the original industrial distribution data of each state and use multi-objective optimization and linear regression to obtain its allocation algorithm.

Finally, we use the mathematical means of the four-dimensional space solution system and the multivariable partial derivative equation to abstract the complex coupled water scheduling relationship into a simple matrix parameter model and use the simulated annealing algorithm to obtain its local optimal solution.

\section{Preliminary} 

\subsection{Assumptions and Justifications}
\begin{itemize}
\item Suppose the lake is a cone whose deepest depth is three times its average depth.

\item 40\% of the water flowing from Glen Canyon Dam goes to Hoover Dam.

\item Water is transported from the dam to the states, and the distance from the dam to the state is regarded as the distance from the dam to the center of the capital.

\item The cost of dams transporting water to the state is only related to distance and elevation difference.
\item The power generation efficiency of the dam at different times is regarded as a fixed value.
\item Water loss only considers state water use, water use for power generation, and water lost by evaporation.
\item The size of a city is determined by its share of GDP in the U.S.
\item The proportion of each industry in a city is determined by its social benefits, economic benefits, and resource consumption.
\end{itemize}

\subsection{Symbols and Notations}
\begin{table}[htb] 
		\centering
        \caption{Notations and Descriptions} 
		\begin{tabular}{ccc} \hline
        \toprule
			Notations&Descriptions&Units\\ 
   \toprule
			$dis_{ij}$&Distance between $dam_i$ and $state_j$&km\\
			$water_{ij}$&Amount of water transported from $dam_i$ to $state_j$&acre-feet\\
            $C_i$&Amount of water that dam $i$ can provide& acre-feet\\
			$R_j$&Amount of water that can be used in $state_j$& acre-feet\\
			$D_j$&Amount of water that can be required in $state_j$& acre-feet\\
            $m,p$&Amount of water in Lake Mead and Lake Powell&acre-feet\\
			$Ele$&Electricity generation& TWh\\
            $\eta$&Electricity generation efficiency&MWh/feet$^2$\\
			$p_i$&Priority in $industry_i$&-\\
            $WF$&Annual water flow of the given river&acre-feet\\
			$t_i$&Proportion of pollutants in sewage&-\\
			$\Delta h$&Dam height difference& m\\

        \bottomrule \hline
	\end{tabular}
\end{table}

\section{Solution to Problem 1}

\subsection{State Water Demand}

As requested by the task, we intend to study Natural resource officials in the US states of Arizona (AZ), California (CA), Wyoming (WY), New Mexico (NM), and Colorado (CO) at two dams and the Colorado River Basin and the water allocation strategies in five states.

The relevant information of these five states was collected from the US government online, including data on population, area, capital, etc., as shown in the table.

\begin{table}[htb] 
    \centering
    \caption{Information about Five States 1} 
    \begin{tabular}{cccc} \hline
    \toprule
    State Name & Abbr. & Total Area& Population\\
\toprule
    Arizona & AZ & $113,998$ sq.mi ($295,254$ $km^2$) & $6,553,255$\\ 
    California & CA & $163,700$ sq.mi ($423,970$ $km^2$) & $38,041,430$\\ 
    Colorado & CO & $104,185$ sq.mi ($269,837$ $km^2$) & $5,187,582$\\ 
    New Mexico & NM & $121,697$ sq.mi ($315,194$ $km^2$) & $2,085,538$\\ 
    Wyoming & WY & $97,818$ sq.mi ($253,348$ $km^2$) & $576,412$\\ 
    \bottomrule \hline
\end{tabular}
\end{table}

\begin{table}[htb] 
    \centering
    \caption{Information about Five States 2} 
    \begin{tabular}{ccccccc} \hline
    \toprule
    State Name&&Abbr. &&Capital && Most Populous City \\ 
\toprule
Arizona &&AZ &&Phoenix &&Phoenix \\ 
California && CA&&Sacramento&&  Los Angeles \\ 
Colorado && CO&&Denver&&  Denver \\ 
New Mexico&& NM&&Santa Fe && Albuquerque \\ 
Wyoming && WY&&Cheyenne && Cheyenne\\ 
\bottomrule \hline
\end{tabular}
\end{table}

We will simplify the distribution model and abstract the distances from the five states to the two dams as the distances from the central cities of the five states to the two dams, and we can get the distances in the table below.

\begin{table}[htb] 
		\centering
        \caption{Dams and States Distance (Unit:km)} 
		\begin{tabular}{ccccc} \hline
        \toprule
			State-Capital&&Glen Canyon dam&&Hoover Dam\\ 
   \toprule
			CA-Sacramento&& $900.63$ && $663.46$\\ 
			AZ-Phoenix&& $391.85$ && $389.35$\\ 
            WY-Cheyenne&& $739.32$ && $1030.90$\\ 
			NM-Santa Fe&& $513.83$ && $792.55$\\
			CO-Denver&& $644.78$ && $956.08$\\
        \bottomrule \hline
	\end{tabular}
\end{table}

According to the relevant data, in 2020, the total water consumption in the United States will be $18$ million acre-feet\cite{c8}.

In 2020, the water supply of the dam is only $5$ million acre-feet. According to the data, the proportion of water obtained from the dam in the five states is

\begin{table}[htb] 
    \centering
    \caption{Water Volume Drawn from Dams in Five States} 
    \begin{tabular}{ccccc} \hline
    \toprule
    State && Abbr. && Proportion of Water\_\text{taken from the dam}\\ 
\toprule
    Arizona && AZ && $17$\% \\ 
    California && CA && $27$\% \\ 
    Wyoming && WY && $6$\% \\ 
    New Mexico && NM && $5$\% \\ 
    Colorado && CO && $23$\% \\ 
    \bottomrule \hline
\end{tabular}
\end{table}

\begin{figure}[h]\label{3111}
\small
\centering
\includegraphics[width=10cm]{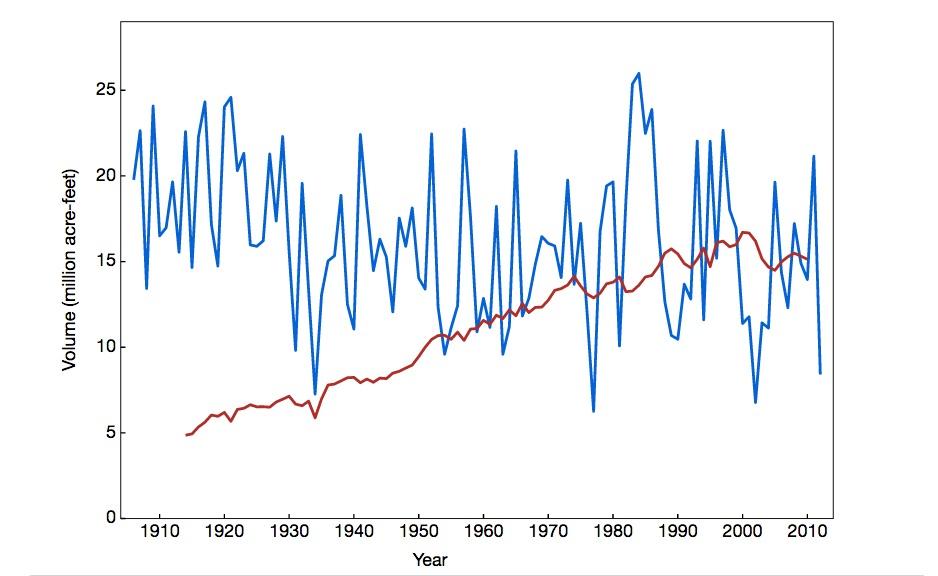}
\caption{Trends in water supply and water use in the Colorado River Basin}
\end{figure}

Therefore, the water demand of each state can be obtained as shown in the table below

\begin{table}[htb] 
    \centering
    \caption{Water Demand in Five States} 
    \begin{tabular}{ccccc} \hline
    \toprule
    State && Abbr. &&Water Demand(MAF)\\ 
\toprule
    Arizona && AZ && $3.06$ \\ 
    California && CA && $4.86$ \\ 
    Wyoming && WY && $1.08$ \\ 
    New Mexico && NM && $0.9$ \\ 
    Colorado && CO && $4.14$ \\ 
    \bottomrule \hline
\end{tabular}
\end{table}

Next, let's look at the calculation formula for hydropower:

The power output of a dam is calculated using the potential energy of the water and can be calculated using the following hydroelectric formula:

\begin{equation} \label{rho}
P = \eta * \rho * g * h * Q
\end{equation}

where P is the power output, measured in Watts, $\eta$ is the efficiency of the turbine, $\rho$ is the density of water, taken as $998$ kg/m³ (you can change it in advanced mode), $g$ is the acceleration of gravity, equal to $9.81$ m/s² (you can change it in advanced mode), $h$ is the head, or the usable fall height, expressed in units of length (meters or feet), $Q$ is the discharge (also called the flow rate), calculated as $Q = A * v$, $A$ is the cross-sectional area of the channel, $v$ is the flow velocity.

Run-of-river installations and tidal power stations take advantage of the kinetic energy of the flow, so the formula is slightly different:

\begin{equation} \label{rho}
P = 0.5 * \eta * \rho * Q * V^2
\end{equation}

The efficiency of the turbine is the ratio of available energy of water to the actual power output of the turbine. It's usually expressed as a percentage. The efficiencies of such turbines can reach up to $59.3$\%, as they're limited by the Betz limit.

\subsection{Dam Water Dispatch}

To meet the water supply needs of the five states, water needs to be diverted from Lake Powell as well as Lake Mead. In the case of meeting the water consumption of residents, we also need to consider the material and manpower consumed by the water resources in the process of transportation. The direct transportation cost of water resources mainly depends on the horizontal distance and height difference between the two locations, and follows the following Principle: For every $1$\% increase in transportation volume, the total cost increases by $0.92$\% and the unit cost decreases by $0.08$\%. Taking into account the life of the transmission pipeline and the loss during transportation, the unit cost of the pipeline fee can be obtained.

Transportation capacity is another major factor that affects the cost of water transportation, which depends on the water scarcity and transportation conditions at the node. Among them, the water shortage degree can be obtained by calculating the water source gap of the target node:

\begin{equation} \label{rho}
C_n = (R_n-D_n)\times p
\end{equation}

where $p$ represents the distance factor, which represents the different modes of transportation.

According to the water transportation data on the official website given by the United States, it can be calculated that $p$ is $0.32$.

According to the literature, we can obtain the formula for calculating the unit cost of water resources:

\begin{equation}
\left\{
             \begin{array}{lr}
             D_c=L_c=5 \;\text{cents}/(\text{m}^3 \cdot100\text{km}) \\
             C_t=d\times D_c + (e_d-e_s)\times L_c\times (1-0.08\%)^{\log_{0.01}^{\frac{C_q}{100}}}  
             \end{array}
\right.
\end{equation}

$e_d$ is the altitude of the destination point, $e_s$ is the altitude of the starting point (unit: 100m)

$c_q$ is the total amount of water transported during water resource scheduling (unit: 100m)

We set the water volume of the Gran Canyon dam to be transported to CA-Sacramento, AZ-Phoenix, WY-Cheyenne, NM-Santa Fe, CO-Denver five states respectively $x_1, x_2, \cdots, x_5$, similar, Hu The water supply from the Buddha Dam to the five states is $x_6,x_7,\cdots,x_{10}$. 

Let the total transportation cost be $W$, $C_1,C_2,\cdots,C_{10}$ are the costs in these $10$ transportation processes, m and p are the water volumes of the two reservoirs, we can establish the relationship between them and the reservoir water level. According to the data, when the water level of Lake Mead is $1067$ ft, the water in the reservoir is only $34$\% of the total capacity. Since the plane area of the reservoir is unchanged, we can set the reservoir capacity and water level as A linear relationship, the total capacity of Lake Mead is $38,296,200,000$ cubic meters, the relationship between water storage and water level can be obtained, and the same is true for Lake Powell's algorithm, the following formula can be obtained

\begin{equation}
\left\{
             \begin{array}{lr}
             m=6.4\times10^8\times M-1.96\times10^{11}  \\
             p=6.78\times10^8\times P-7.01\times10^{11} 
             \end{array}
\right.
\end{equation}

The water supply needs of the residents of the five states were previously derived as S1...S5. We can calculate the altitude of the water level of the reservoir as ed, and the altitudes of the five states are as follows:

\begin{table}[htb] 
    \centering
    \caption{Altitudes of the Five States} 
    \begin{tabular}{cccccc} \hline
    \toprule
    &State - Capital City&&Altitude(m)&\\ 
\toprule
&CA - Sacramento && $9.13$&\\ 
&AZ - Phoenix  && $331$&\\ 
&WY - Cheyenne  && $1848$&\\ 
&NM - Santa Fe && $2194$&\\ 
&CO - Denver  && $1596$&\\ 
\bottomrule \hline
\end{tabular}
\end{table}

We use the following constraints to bind w to find the minimum value of w: the total amount of water supply should not be greater than the capacity of the reservoir, and secondly, we need to meet the water demand of the households, and list the equations

\begin{equation}\label{WC}
\left\{
             \begin{array}{lr}
             W=\sum_{i=1}^{10}C_i \\
              C_t=d_i\times D_c + (e_d-e_s)\times L_c\times (1-0.08\%)^{\log_{0.01}^{\frac{C_q}{100}}} \\
              \sum_{i=1}^5 x_i\le p,\;\;\sum_{i=6}^{10} x_i\le m \\
              x_k+x_{k+5}=s_k,\;\;k=1,2,\cdots,5
             \end{array}
\right.
\end{equation}

We can solve equation(\ref{WC}) by MATLAB.

Under the constraints we can derive the allocation of water volume to each state for each reservoir at a given water level.

Therefore, inputting different M and P can get the following water distribution.


\begin{table}[htb] 
    \centering
    \caption{Altitudes of the Five States} 
    \begin{tabular}{cccccccccccc} \hline
    \toprule
    M(m)&P(m)&x\_1(MAF)&x\_2&x\_3&x\_4&x\_5&x\_6&x\_7&x\_8&x\_9&x\_{10}\\ 
\toprule
    1067.18& 3529.37 &0 & 0 &0.9&0.75&3.45&2.55&4.05&0&0&0\\ 
    1020 & 3800 &0& 0&0&0&2&2.55&4.05&0.9&0.75&1.45\\
    1180 & 3500 &1.78& 4.05&0.9&0.75&3.45&0.77&0&0&0&0\\
\bottomrule \hline
\end{tabular}
\end{table}

According to the above data, we can get that, under the actual conditions of the current reservoir, Lake Powell should pump $0.9$ MAF (million acre-feet) of water to Wyoming, $0.75$ MAF to New Mexico, and $3.45$ MAF each year. For Colorado, a total of $5.1$ MAF of water was withdrawn. Lake Mead should pump $2.55$ MAF of water annually to New Mexico and $4.05$ MAF to Arizona, for a total of $6.6$ MAF of water.

The Colorado originates in the Rocky mountains and traverses seven US states, watering cities and farmland, before reaching Mexico, where it is supposed to flow onwards to the Sea of Cortez. Instead, the river is dammed at the US-Mexico border, and on the other side the river channel is empty. The United States built a series of dams and canals on the Colorado River to trap all the river water, so the river has not flowed to the sea since 1998\cite{c6}. 

\begin{table}[htb] 
    \centering
    \caption{Average Urban Water Use in CRD (Mexicali)} 
    \begin{tabular}{ccccc} \hline
    \toprule
    &Use cases&Consumption($1,000 m^3/year$)&Share(\%)& \\
\toprule
&Residential&$57,125$&$73.3$& \\
&Commercial&$7,197$&$9.2$& \\
&Industrial&$6,352$&$8.2$&\\
&Other&$7,234$&$9.3$&\\
&Total&$77,908$&$100.0$&\\
\bottomrule \hline
\end{tabular}
\end{table}

As mentioned above, our model is analyzing the agricultural, residential, and industrial water use in Mexico, collecting relevant data, and combining it with the distribution data of five states in the United States to further improve the distribution model and consider the interests of Mexico.

\subsection{Calculation of Water for Power Generation}

According to the principle of hydropower generation, we can conclude that dam power generation is essentially the conversion of the gravitational potential energy of water into electrical energy. 

From the physical potential energy formula(\ref{E}),
\begin{equation} \label{E}
\Delta E=E_{p1}-E_{p2}=mgh_1-mgh_2=mg\Delta h
\end{equation}
we guess that the power generation is not only linearly related to the water flow, but also linearly related to the drop of the dam.

Taking Lake Powell as an example, the data query shows that its current water level is $3528.72\;\text{feet}$, and the maximum water level is $3700\;\text{feet }$, the maximum depth of Lake Powell is $558\;\text{feet}$, then it can be inferred that its bottom is $3142\;\text{feet}$, Its current lake depth $\Delta h=3528.72-3142=386.72\;\text{feet}$.

The following verify that there is a linear relationship between power generation and water flow and dam drop:

By analyzing the data, it can be found that between $2017-2021$ (in order to make the data representative and fit the actual situation, the data of the last 5 years was selected for analysis), the average water level of Lake Powell is $3605\ ;\text{feet}$, so select the data between $3600-3610\;\text{feet}$ to analyze the correlation between water flow and power generation.

\begin{figure}[h]\label{31}
\small
\centering
\includegraphics[width=10cm]{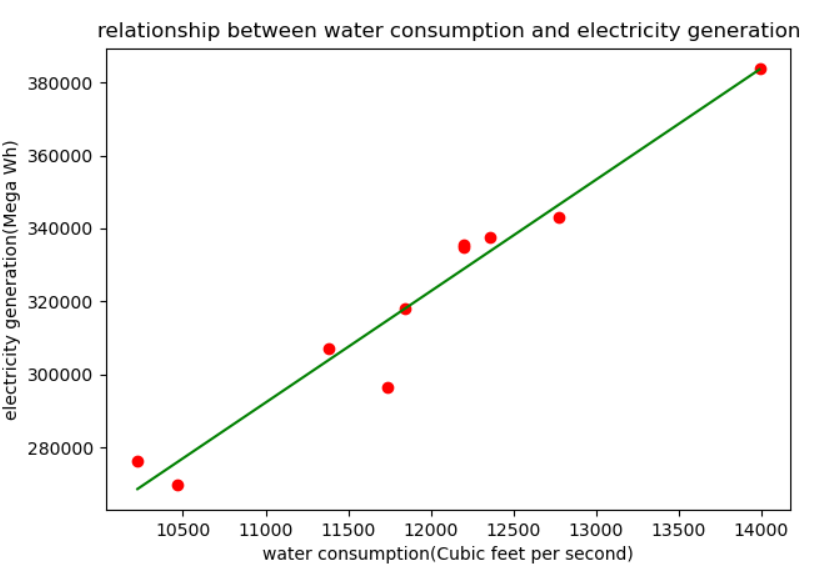}
\caption{Relationship between water consumption and electricity generation}
\end{figure}

It is not difficult to find a highly linear correlation between water flow and power generation

The next step is to verify the relationship between power generation and water potential drop.

We only need to verify the relationship between the ratio of power generation and water flow $\frac{electricity\;generation}{water \;consumption}$ and the water potential drop. According to the data, the relationship can be obtained as follows:

\begin{figure}[h]\label{32}
\small
\centering
\includegraphics[width=10cm]{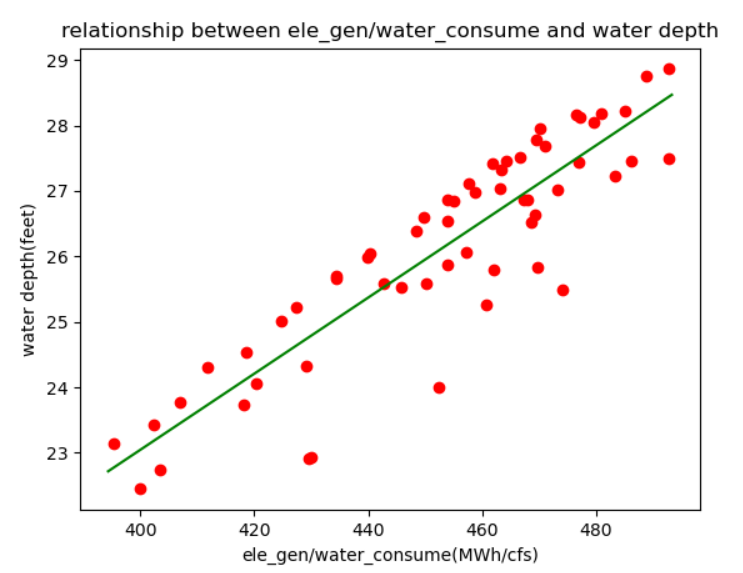}
\caption{Relationship between ele\_gen/water\_consume and water depth}
\end{figure}

It is not difficult to find that there is also a good linear correlation between the ratio of power generation and water flow and the water potential drop.

From this, we can get the relationship between the power generation $Ele$ and the water potential drop $\Delta h$ and the water flow $WF(water\;flow)$ as:

\begin{equation} 
Ele=\eta\times\Delta  h \times WF
\end{equation}

Among them, $\eta$ is the mechanical efficiency of power generation. Substitute the data of the past $5$ years to calculate $\eta$ and take the average value to get

\begin{equation}
\eta = 5.75 \times 10^{-2}  (\text{MWh}/\text{feet}^2)
\end{equation}

(The standard deviation of $\eta$ in 60 samples in $5$ years is $1.9\times 10^{-3}$, which also shows that the bilinear model has a good fit)

Finally, combining the information and data on the Internet, we can obtain the water consumption for power generation in Lake Mead and Lake Powell as follows:

\begin{figure}[h]\label{33}
\small
\centering
\includegraphics[width=10cm]{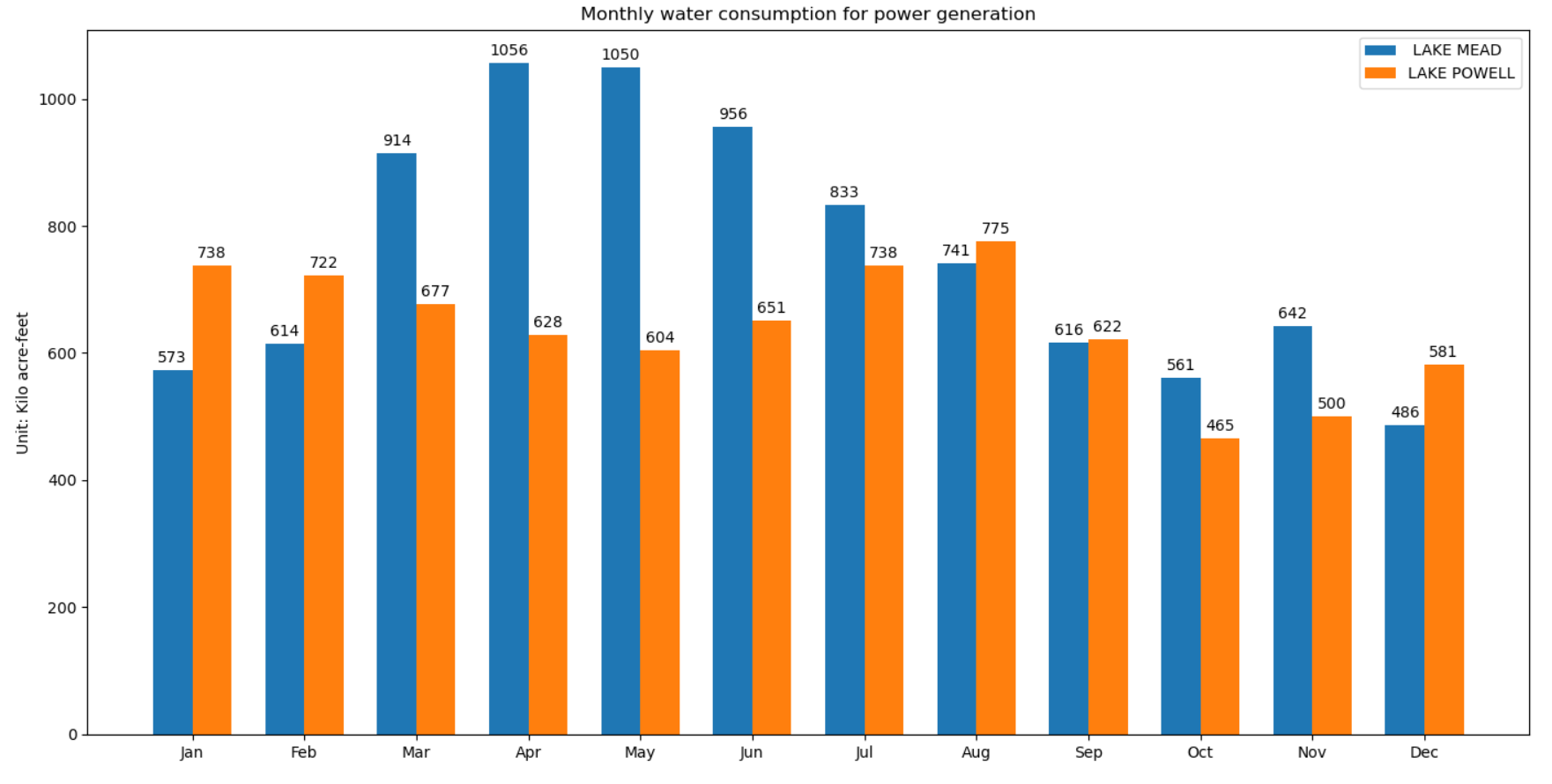}
\caption{Monthly water consumption for power generation}
\end{figure}

\subsection{Evaporation of Water}

According to the collected data, it can be seen that the monthly evaporation in the last ten years has not changed much, so we can use the average value of the data in the past ten years as the predicted value of the future evaporation.

\begin{table}[htb] 
    \centering
    \caption{Average Evaporation(2010-2020)(Unit:AF)} 
    \begin{tabular}{cccccccccc} \hline
    \toprule
   Time&AE of Lake Powell&AE of Lake Mead\\ 
\toprule
    1-Jan & $7847.9$&$7156$\\ 
  1-Feb & $8272$&$7543$\\
   1-Mar & $13890.2$&$12666$\\
  1-Apr  & $21932.6$&$19999$\\ 
 1-May & $26450.3$&$24119 $\\ 
   1-Jun & $44017.9$& $ 40138$\\
  1-July & $54031.5$&$49269$\\ 
  1-Aug  & $52619.4$&$47981$\\
1-Sep  & $47521.8$& $43333$\\ 
 1-Oct  & $32587.4$& $29715$\\ 
  1-Nov  & $31048.8$&$28312$\\ 
1-Dec  & $24297.8$&$22156$\\ 
    \bottomrule \hline
\end{tabular}
\end{table}

 Evaporation increased sharply in January and July, reached its peak in July, and then decreased year by year as the years went on. The U.S. Southwest has experienced its worst drought in nearly two decades since the summer, leading to a continued decline in water resources in the Colorado River Basin, leaving the country's two largest reservoirs at record low levels. The US Oceanic and Atmospheric Administration believes that the weak southwest monsoon in the summer of 2020 has resulted in less precipitation in the southwest of the United States since 2020. In addition, the continuous high temperature weather in this area in the summer of 2021 will intensify evaporation, especially in the summer. volume increased sharply. The overall evaporation remains above $10,000$ AF, which is a relatively high level.
 
 Let's do the evaporation analysis for Lake Mead (similar to above).
 

\subsection{Lake Maintenance Time}
According to the U.S. Bureau of Reclamation (USBR), on February 14, 2022, Lake Powell's water volume was only $25\%$ of its full capacity, and its current water level was $3529.37\text{feet}$, with a total capacity of $3529.37\text{feet}$ is $6199\;\text{KAF(kilo acre-feet)}$, while for Lake Mead this figure is $8977\;\text{KAF}$.

For a water storage system composed of lakes and dams, the main sources of water are surface water inflow, upstream water inflow of lakes, rainfall, and melting of glaciers and snow, and the main sources of its consumption are evaporation loss, water for power generation, and water for state supply into the downstream reaches, etc\cite{c7}.

\begin{figure}[h]\label{34}
\small
\centering
\includegraphics[width=10cm]{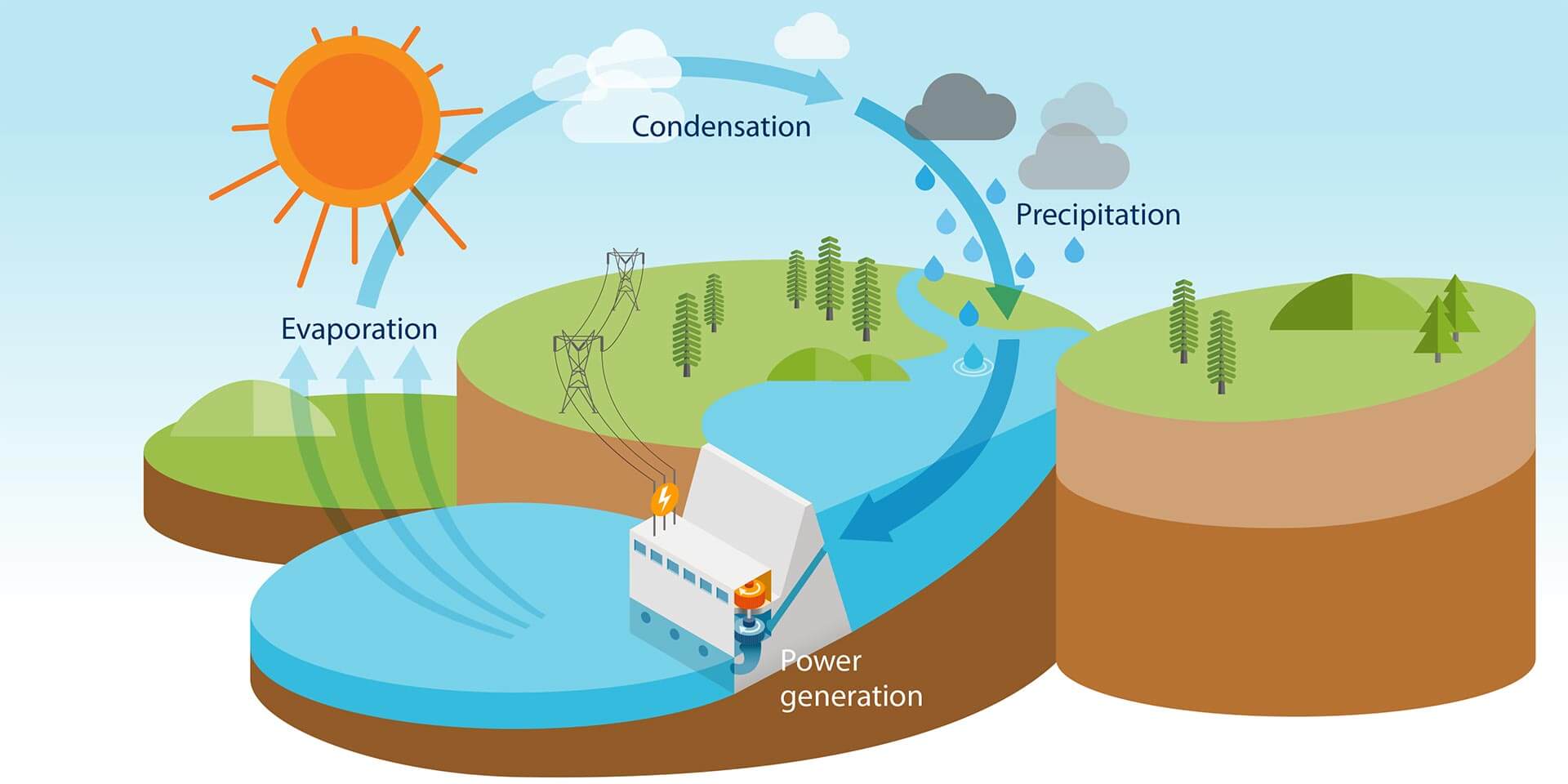}
\caption{Simple diagram of water cycle}
\end{figure}

From the above article, we discuss the three main paths of water consumption: evaporative losses, state water use, and water use for power generation.

According to the data of several states in 2020 obtained by our search, Lake Powell and Lake Mead supply a total of $11.7$ MAF of water to these five states, which is an average of $0.975$ MAF per month, and the total monthly total of Lake Powell and Lake Mead The evaporation and power generation water consumption are shown in the following figure:

\begin{figure}[h]\label{35}
\small
\centering
\includegraphics[width=12cm]{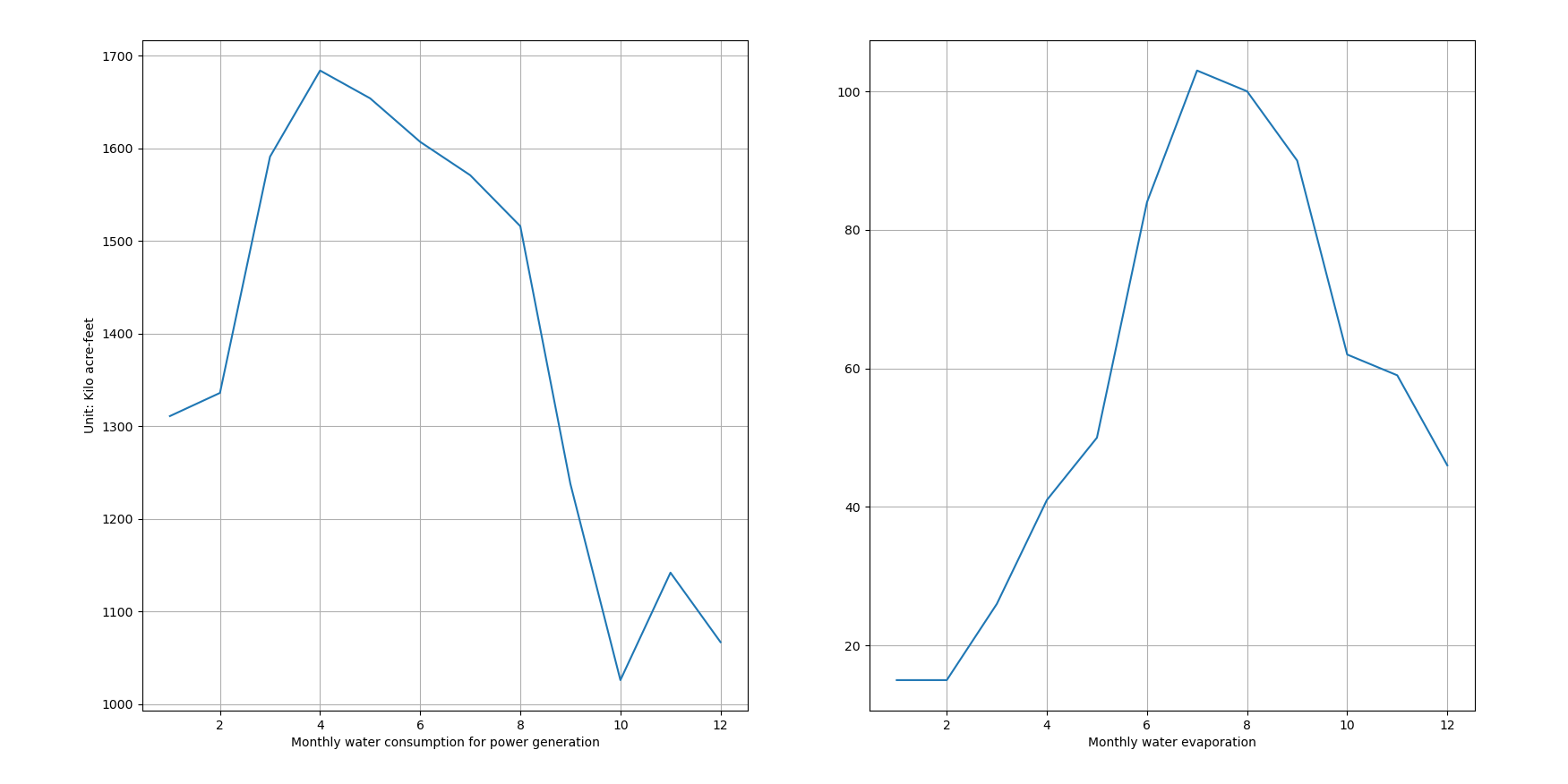}
\caption{Monthly water consumption for power generation and monthly water evaporation}
\end{figure}

According to the previous state water use, power generation water, and evaporation water, we conclude that it will be used up after $5.76$ months, that is, after about $173$ days.

\subsection{Additional Water Supply}
According to the data in the previous question, if we want to maintain the supply of water, that is, to make up all the consumption, then the amount of water we need every year is the sum of the above three, that is, the sum of state water, power generation water and evaporation water, its value $29.13$ MAF

\begin{figure}[h]\label{36}
\small
\centering
\includegraphics[width=10cm]{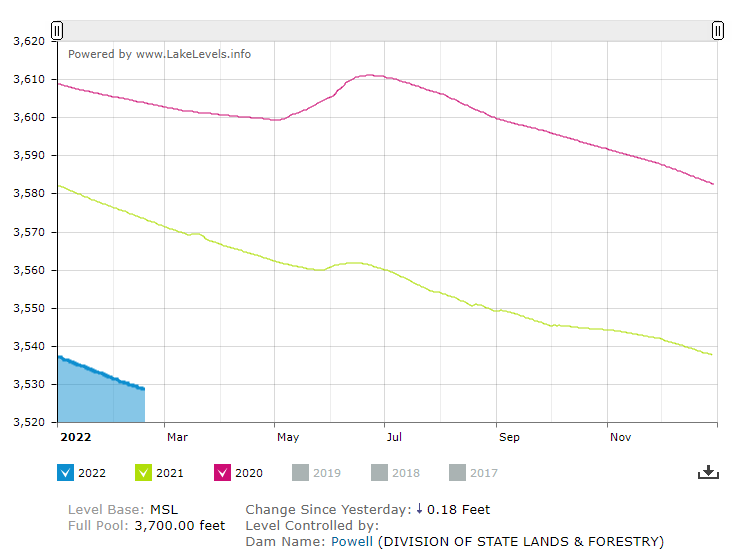}
\caption{Reduced water volume of Lake Powell(KAF}
\end{figure}

If the precipitation and inflow water are considered here, it can be seen from the data query that Lake Powell will reduce the water volume by $330$ KAF in 2021, and Lake Mead will reduce the water volume by $1410$ KAF in 2021, which means that a total of $1.74$ MAF needs to be added. water volume.

\section{Solution to Problem 2}
\subsection{Model Building}

According to the requirements of the question, we need to allocate the water resources for general (agricultural, industrial, residential) use and for the production of electricity under the given the water resources.

For the allocation of water resources we need to achieve the best overall efficiency. The first is the construction of the objective function, which should include the maximum economic benefit, the maximum social benefit and the maximum environmental benefit.

\subsubsection{Economic Efficiency Goal}
The economic benefits mainly consider the economic benefits generated by the water and electricity consumption of agriculture and industry. We set w as the total economic benefit and the water supply for agriculture is $c_1$, the power supply for agriculture is $c_2$, the water supply for industry is $c_3$, and the power supply for industry is c4. So there are:

\begin{equation}
\max W = \sum_{i=1}^4 e_iC_i
\end{equation}

In the formula $e_i$ is unknown, and we predict the parameter values of $e_i$ by building a gray model.
First of all, we found the data:
According to the U.S. Bureau of Economic Analysis (BEA), the total U.S. GDP in 2021 is $22.99$ trillion dollars, and the GDP of the five states given in the question are

\begin{table}[htb] 
    \centering
    \caption{State's 2021 GDP(Unit: million of US dollars)} 
    \begin{tabular}{cc} \hline
    \toprule
    Unit: million of US dollars&state's 2021 GDP\\ 
\toprule
California (CA) & $3,353,473 $\\ 
Colorado (CO)  & $425,595 $\\ 
Arizona (AZ)  & $409,577 $\\ 
New Mexico (NM)  & $110,696 $\\ 
Wyoming (WY) & $42,536 $\\  
    \bottomrule \hline
\end{tabular}
\end{table}

Also according to the Central Intelligence Agency, we can know that the shares of agriculture, industry, and services in the United States are shown below.

\begin{table}[htb] 
    \centering
    \caption{Share of US industries} 
    \begin{tabular}{ccc} \hline
    \toprule
    Agriculture&Industry&Services\\ 
\toprule
    0.9\% & 18.9\% & 80.2\% \\ 
    \bottomrule \hline
\end{tabular}
\end{table}

This gives us a total agricultural output of $39,077$ million dollars/year and a total industrial output of $820,614$ million dollars/year for the five states.

Based on data of Year 2019, the following graph shows the consumption of electricity in the United States by channel

\begin{table}[htb] 
    \centering
    \caption{2021 U.S. Electricity Consumption by Industry} 
    \begin{tabular}{ccccccc} \hline
    \toprule
    &Industry&Transport&Residential&Commercial&Other&Total \\ 
\toprule
    Unit: TWh & $9566$ & $420$ & $6072$ & $4849$ & $1940$ & $22847$ \\ 
    \bottomrule \hline
\end{tabular}
\end{table}

The calculated share of industrial electricity in total electricity consumption in the United States is about $42.5$\%, and the share of residential electricity in total electricity consumption is about $26.6$\%. And according to the U.S. DEPARTMENT OF AGRICULTURE (USDA), U.S. agricultural electricity use is about $1.74$\% of total electricity useCombining the data in the first question, we can see that these five U.S. states account for approximately 18.96\%.

From this, we can estimate the annual electricity consumption of these five states for agriculture, industry, and residential use as

\begin{table}[htb] 
    \centering
    \caption{Electricity consumption by industry in five states} 
    \begin{tabular}{cccc} \hline
    \toprule
    &Industry&Residential&Agriculture \\ 
\toprule
    Unit: TWh & $1813$ & $1151 $ & $76$ \\ 
    \bottomrule \hline
\end{tabular}
\end{table}

We can get from the United States Geological Survey (USGS), the data for agricultural, industrial, and residential water use in these five states are shown below\cite{c9}.

\begin{table}[htb] 
    \centering
    \caption{Water Consumption by Industry in Five States} 
    \begin{tabular}{cccc} \hline
    \toprule
    Unit: M gallons/day&Agriculture&Industry&Residential \\ 
\toprule
    California (CA) & $19,000$ & $399$ & $127$ \\ 
    Colorado (CO) & $9,000 $ & $84.1$ & $35.4 $ \\ 
    Arizona (AZ) & $4,530$ & $6.12 $ & $24.0 $ \\ 
    New Mexico (NM) & $2,370$ & $24.6$ & $3.40 $ \\ 
    Wyoming (WY)  & $7,790$ & $8.04 $ & $8.93 $ \\ 
    \bottomrule \hline
\end{tabular}
\end{table}

The total agricultural, industrial, and residential water use and share in these five states can be found as follows

\begin{table}[htb] 
    \centering
    \caption{Total Water Consumption by Industry } 
    \begin{tabular}{cccc} \hline
    \toprule
    &Agriculture&Industry&Residential \\ 
\toprule
    Unit: M gallons/day & $42690$ & $522$ & $198$ \\ 
    percentage(\%) & $98.3$ & $1.2$ & $0.5 $ \\ 
    \bottomrule \hline
\end{tabular}
\end{table}

Based on the relationship between agricultural and industrial GDP and water supply and electricity supply in each state in 2021, we can thus obtain the coefficient of return to water supply at different water supply levels $e_1$.

For the prediction of $e_1$ we use a gray model for prediction, and we can obtain the relationship between water supply and the coefficient of return of water supply as follows.

We can obtain the series from the above data as $e_1^{(0)}$, which are $e_{11}$, $e_{12}$, $e_{13}$, $e_{14}$, $e_{15}$. These numbers are the water supply revenue coefficients for the five states, and we accumulate their values as a new series, that is
\begin{equation}
e_1^{(1)}(k)=\sum_{i=1}^ke_1^{(0)}(k)
\end{equation}

After that we also obtain the background value from $e_1^{(0)}$,that is 
\begin{equation}
\lambda^{(1)}(k+1)=\frac12[e_1^{(1)}(k+1)+e_1^{(1)}(k)]
\end{equation}
so there are
\begin{equation}
\left\{
             \begin{array}{lr}
			e_1^{(0)}(k)+\varepsilon\lambda^{(1)}(k)=\mu[\lambda^{(1)}(k)]^2 \\
			\frac{de_1^{(1)}}{ds}+\varepsilon e_1^{(1)}=\mu(e_1^{(1)})^2
             \end{array}
\right.
\end{equation}

Where $s$ is the water supply
$\varepsilon$ and $\mu$ are development coefficients and grey predictors respectively
The least square estimate of $\varepsilon$ is $\hat{\varepsilon}$ 
\begin{equation}
\hat{\varepsilon}=[\varepsilon,\mu]^T=(H_v^TH_v)^{-1}H_v^TH_v
\end{equation}

\begin{equation}
H_v=\begin{bmatrix}
  -\lambda^{(1)}(2)& -\lambda^{(1)}(2)^2\\
  -\lambda^{(1)}(3)& -\lambda^{(1)}(3)^2 \\
  -\lambda^{(1)}(4)& -\lambda^{(1)}(4)^2 \\
  -\lambda^{(1)}(5)& -\lambda^{(1)}(5)^2
\end{bmatrix}
\end{equation}

$\varepsilon$ and $\mu$ are calculated according to the formula above

The grey model time response in the initial state can be combined with $\varepsilon,\mu$ and the formula above, it can be solved that
\begin{equation}
\hat{e_1}^{(1)}(k+1)=\frac{\varepsilon e^{(1)}(1)}{\mu e^{(1)}(1)+[\varepsilon-\mu e^{(1)}(1)]k^{sk}}\;\;k=2,3,4,5
\end{equation}

Finally, the $\hat{e_1}^{(1)}$ sequence can be reduced to $\hat{e_1}^{(0)}$, from which the water supply income coefficient under different water supply can be obtained

We can obtain the trend graph of the corresponding benefit coefficients

The other revenue coefficients can also be found in the corresponding way. The revenue coefficient of agricultural power supply has been at a low level, and the trend of the other revenue coefficients is similar to the revenue coefficient of agricultural water supply, which basically increases rapidly with the amount of water supply first, and then increases slowly, and finally the revenue coefficient basically remains unchanged when the amount of water supply reaches a certain level.

\subsubsection{Social Benefit Goal}

The social benefit objective is also the minimum water and electricity shortage for agriculture, industry and housing. We set $f_1$as the water shortage and $f_2$ as the electricity shortage, then there are:

\begin{equation}
\min f_1=o_1-c_1+o_3-c_3+o_5-c_5
\end{equation}
$o_1,o_3,o_5$ are the water demand for agriculture, industry, and residential
\begin{equation}
\min f_2=o_2-c_2+o_4-c_4+o_6-c_6
\end{equation}
$o_2,o_4,o_6$ are the electricity demand for agriculture, industry, and residential
$c_5,c_6$ are the water and electricity consumption of the house

We can know $O_i$ ($i=1,2,3,4,5,6$) according the data.
According to a survey conducted by the U.S. Bureau of Official Statistics\cite{c}10, for agriculture, industry, and people's livelihood, the satisfaction of water demand and electricity demand are 

\begin{table}[htb] \label{t18}
    \centering
    \caption{Satisfaction with Water and Electricity Supply by Industry} 
    \begin{tabular}{cccc} \hline
    \toprule
    &Agriculture&Industry&Residential \\ 
\toprule
    Water demand satisfaction & 96\% & 73\%  & 84\%  \\ 
    Electricity demand satisfaction & 56\%  & 88\% & 97\%  \\ 
    \bottomrule \hline
\end{tabular}
\end{table}

So we can get the actual demand of water and electricity for agriculture, industry and residents as:

\begin{table}[htb] 
    \centering
    \caption{Actual Water and Electricity Demand by Industry} 
    \begin{tabular}{cccc} \hline
    \toprule
    &Agriculture&Industry&Residential \\ 
\toprule
    Unit: TWh & 135 & 2483  & 1370   \\ 
    Unit: M gallons/day & 44468  & 715 & 235 \\ 
    \bottomrule \hline
\end{tabular}
\end{table}

\subsubsection{Environmental Efficiency Goal}

For the environmental benefits, we mainly use the total amount of pollutants in wastewater discharged from agriculture, industry and residential as an indicator. We let the total amount of pollutants be $f_3$ and so there are:

\begin{equation}
\min f_3 = \sum_{i=1}^5c_i\cdot p_i\cdot t_i
\end{equation}
$p_i$ is the ratio of sewage to water, $t_i$ is the proportion of pollutants in sewage

From the data we can calculate $p_i$ and $t_i$.

According to the data of the American Water Quality Survey Network, the volume and concentration of wastewater generated by the three industries of agriculture, industry, and livelihood are as follows.

\begin{table}[htb] 
    \centering
    \caption{Production of Pollutants by Industry} 
    \begin{tabular}{cccc} \hline
    \toprule
    &Agriculture&Industry&Residential \\ 
\toprule
    polluted water quantity(M gallons/day) & 8439 & 690  & 117  \\ 
    polluted water concentration(\%) & 4\%  & 56\% & 32\% \\ 
    \bottomrule \hline
\end{tabular}
\end{table}

When we give the total amount of water resources as $i$, the relationship between water flow and power generation is given by the conclusion of the first question as Equation \ref{26}

\begin{equation}\label{26}
Ele= \eta \times \Delta h \times WF,\;\eta=5.75\times10^{-2} MWh/feet^2
\end{equation}

from which we can obtain the constraint that

\begin{equation}
C_1+C_2+C_3+\frac{C_4+C_5+C_6}{\eta\cdot \Delta h}\le i
\end{equation}
where $\Delta h$ can be taken as the annual average water potential drop of Lake Mead.

\subsection{Solution of the Model and Conclusion}

We use the multi-objective ant colony genetic algorithm to solve the above model, and use this model to solve the competition problem of water and electricity that Lake Mead and Lake Powell should allocate to the five states, we can get the following conclusion, input the total water of the two reservoirs, we can get the annual water supply of Lake Powell and Lake Mead to the five states is $10.53$ MAF for agriculture, $0.702$ MAF for industry, and $16.2$ MAF for power generation\cite{c12}. The amount of water supplied for residential use is $0.468$ MAF, and the amount of water supplied for power generation is $16.2$ MAF.

\section{Solution to Problem 3 \& 4}
\subsection{Problem Analysis}
This problem can be analyzed in two levels. First, according to the distance from the reservoir to the state, calculate how much water the reservoir should allocate to each state. It should be clear here that the water demand of the state and the degree of satisfaction are not linear. For example, when When a state is supplied with $90$\% of its water needs, the effect on the state may not be much because crops will not stop growing because of a little less water, and if a state is supplied with only $50$\% of its water needs, then The operation of the state will be greatly affected and even the state will be shut down due to lack of water.

After the water is transported from the reservoir to the state, the state will allocate the corresponding water resources according to the characteristics of its development. For example, for Arizona, its industry only accounts for a very small part, and agriculture accounts for a large proportion, then We can try to meet its core industrial needs while reducing the amount of its agricultural water supply.

\subsection{Model Establishment}

Regarding the state water demand, we refer to the first-order constant coefficient differential equation in the ordinary differential equation, and denote the water supply and satisfaction of the five cities as $x_i,f_i(x_i),i=1,2,\cdots, 5$, and record its minimum and maximum water consumption as $X_{i\;min}, X_{i\;max}$ to obtain the following reference equation.

\begin{equation} 
\left\{  
             \begin{array}{lr}  
			\frac{df_i(x_i)}{dx_i}=k_ix_i,\;\;x_i\in[X_{i\;min},X_{i\;max}] \\
			f_i(X_{i\;min}) = 0,f_i(X_{i\;max}) = 1, i=1,2,\cdots,5
             \end{array}  
\right.
\end{equation}

From this, we can solve the satisfaction coefficient of state water use $k_1,k_2,\cdots,k_5$

Taking the factor of distance into consideration, it can be seen from the solution process of the first question that there is a close linear relationship between the cost of transporting water and the distance between the reservoir and the city. Our purpose is to optimize the city’s satisfaction as much as possible. To reduce costs as much as possible, we can formulate the following equations, where $p$ is the proportionality constant

\begin{equation}
\left\{  
             \begin{array}{lr}  
			cost_{ij} = water_{ij}\times dis_{ij}\times p \\ 
			COST = \sum_i \sum_j cost_{ij} \\
			SATISFY = \sum_if_i(x_i)\\
			\min COST, \min SATISFY
             \end{array}  
\right.
\end{equation}

From this we can solve for the amount of water $water_{ij}$ that the $j$ lake will supply to the $i$ lake

After water is allocated to cities, we then allocate water to various industries according to the characteristics of the city. Our rule is to allocate water to cities according to the proportion of industries, but the proportion is related to its scale. We make the agricultural water use of each city, Industrial water, residential water, and power generation water are listed as an array, that is, $a=[agriculture,industry,resident,electricity]$, then the city's priority coefficient $p_i,i=1,2,3,4$ is calculated as follows:

\begin{equation}
\left\{
             \begin{array}{lr}
             tot = \sum_ia_i,i=1,2,3,4  \\
             priority_i = \sqrt{\frac{a_i}{tot}}\\
             p_i = \frac{priority_i}{\sum priority}  
             \end{array}
\right.
\end{equation}

\subsubsection{Calculate and Solve}

Since the above formula involves multi-objective optimization, it is not an ordinary linear programming, and its global optimal solution cannot be obtained. Therefore, it is considered to use simulated annealing algorithm to discretize the numerical value to obtain its local optimal solution\cite{c11}.

In the case of a $20$\% reduction in the sum of state water use and power generation water use, the local solution is as follows:


\begin{figure}[h]\label{211}
\small
\centering
\includegraphics[width=10cm]{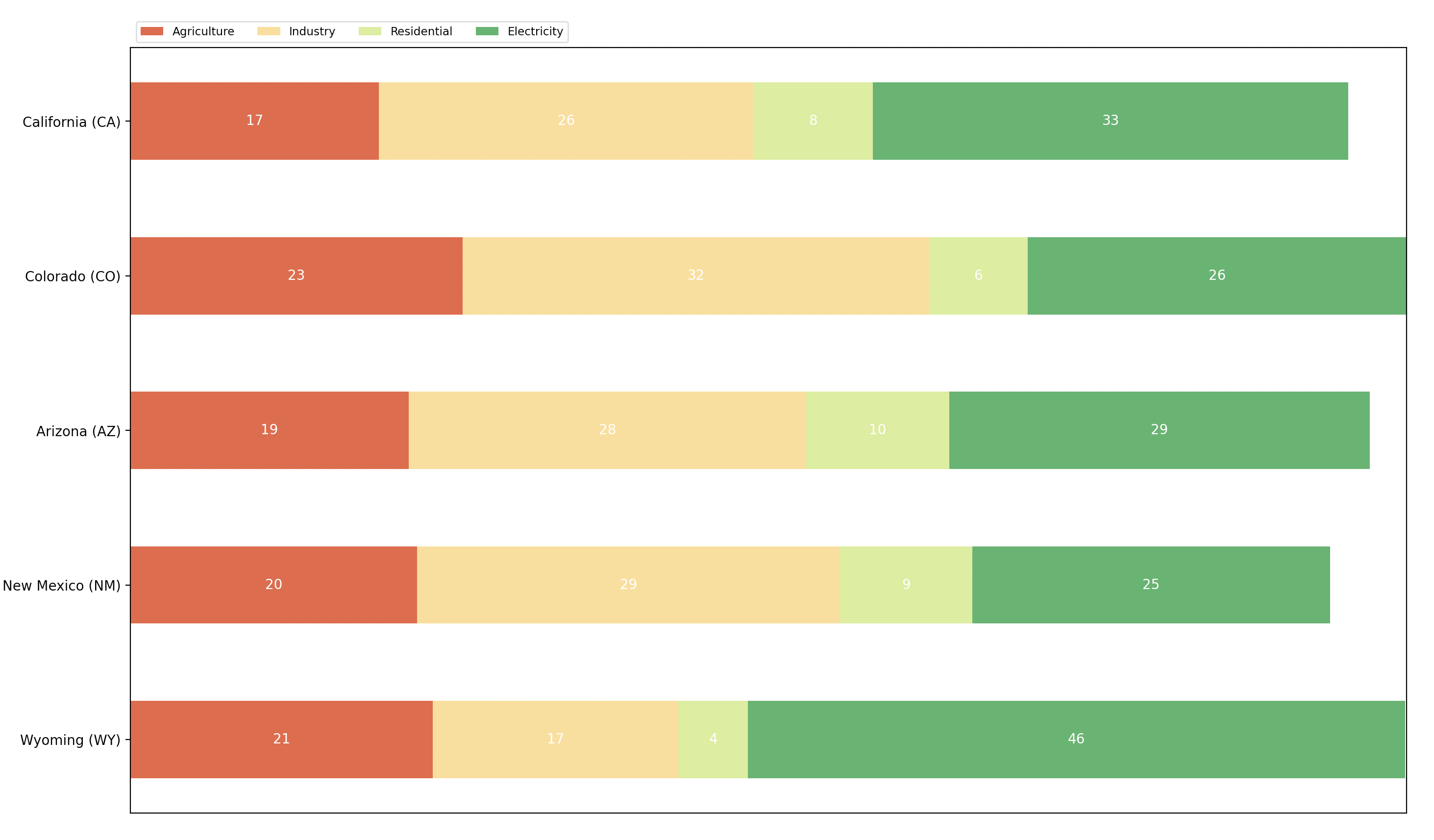}
\caption{Reductions by industry in five states}
\end{figure}

\subsection{Sensitivity Analysis Based on Problem 4}
In view of the three situations proposed in Problem 4, it is essentially to let us solve the final distribution results after the four indicators of industrial water, agricultural water, residential water, and power generation water change. We regard the model as a four-dimensional space. The amount of water allocated to each state is determined by the four dimensions of industrial water, agricultural water, residential water, and water demand for power generation, namely

\begin{equation}
MODEL\_RESULT = f(agriculture,industry,resident,electricity)
\end{equation}

And each indicator change will bring different impacts. For example, if there is less water used for power generation\cite{c14}, then we can use more water for agricultural water, industrial water, and residential water, and this part of water does not necessarily need to be fully allocated to Other uses can be allocated in the order of priority requested before. From this idea, we can list the following formulas

\begin{equation}
\Delta a_j(supply)= \frac{\partial f}{\partial a_i}\times\Delta a_i(demand)\times priority_j\times k,\;\;j = \{1,2,3,4\}-{i}
\end{equation}


For reasons of space, only the data of California is listed below.
According to the calculation of the data, it can be obtained that for case 1, if the population increases by $5$\%, the agricultural scale decreases by $10$\%, and the industrial scale increases by $5$\%, the residential water consumption increases by $4.57$\%, and the agricultural scale increases by $4.57$\%. Water use is reduced by $8.91$\%, industrial water use is increased by $6.70$\%, and power generation water is reduced by $0.98$\%. For case 2, the proportion of renewable energy technologies has increased, that is, power generation water use is reduced. If it is reduced by $10$\%, then residential water use will increase by $3.55$\%, and agricultural water use will increase by $3.55$\%. Water use increases by $0.81$\%, and industrial water increases by $3.4$\%. For case 3, if its demand for hydropower is reduced by $10$\%, then its agricultural water consumption is reduced by $11.86$\%, industrial water consumption is reduced by $8.98$\%, residential water consumption is reduced by $10.06$\%, and power generation water consumption is reduced by $7.57$\%.

\section{Model Evaluation}
\subsection{Strengths}
Our main model solves the most basic problems firstly, and then makes some improvements to the model according to different situations, so that the sub-model can have better practicability. Our models can also get a relatively correct result when it does not need too much data. If the accuracy and quantity of the data are sufficient, the accuracy of the model can also be improved. Our model uses raw data so that the model can reflect the real world situation. The data processing method used by our model helps to avoid some extreme data, which can make the model more robust and universal.

\subsection{Weaknesses}
We mainly use linear data in data prediction, and do not use some algorithms such as neural networks, so the data obtained may have certain errors with the real results.

\subsection{Possible Improvements}
Use higher-precision data to add more parameters to the model, take more account of the real-world economic and political landscape when processing data, and discover the correlations behind the data.


\end{document}